\documentclass[12pt]{article}
\usepackage{mathptmx}

\usepackage{graphicx}
\usepackage{bm}
\usepackage{amssymb}
\makeatletter
\author{Robert P. C. de Marrais \footnote{Email address:  rdemarrais@alum.mit.edu} \\ \noindent
\emph{Thothic Technology Partners, P.O.Box 3083, Plymouth MA 02361}}

\title{Placeholder Substructures I:~~  The Road from NKS to Scale-Free Networks
is Paved with Zero Divisors}

\begin{document}
\maketitle
\makeatother

\begin{abstract}
Zero-divisors (ZDs) derived by Cayley-Dickson Process (CDP) from
N-dimensional hypercomplex numbers ($N$ a power of $2$, and at least
$4$) can represent singularities and, as $N \rightarrow \infty$,
fractals -- and thereby, scale-free networks.  Any integer $> 8$ and
not a power of $2$ generates a meta-fractal or \textit{Sky} when it
is interpreted as the \textit{strut constant} (S) of an ensemble of
octahedral vertex figures called \textit{Box-Kites} (the fundamental
ZD building blocks). Remarkably simple bit-manipulation rules or
\textit{recipes} provide tools for transforming one fractal genus
into others within the context of Wolfram's Class 4 complexity.
\end{abstract}

\section{Introduction:  A Sky for Box-Kites}
Within a generation of Imaginary Numbers' taming, William Rowan
Hamilton etched his famous $i^{2}= j^{2} = k^{2} = ijk = -1$
equation on the bridge he was crossing, thereby giving us
Quaternions (and, eventually, the tools of vector calculus derived
from their workings). Within months, Arthur Cayley and John Graves
independently generalized his 4-D (noncommutative) Complex
quantities to the 8-D (nonassociative) Octonions.  But before the
century was out, attempts to extend hypercomplex number systems
indefinitely, by the CDP dimension-doubling algorithm, met with an
impediment: in 1896, Adolf Hurwitz proved that the next redoubling
(to the 16-D Sedenions) would continue the pattern of giving up
something familiar as the price of admission; only this time, the
very notion of a norm (and hence, of division algebra itself) would
be forfeit. [1] The necessity of zero-divisors (non-null numbers
with zero product) in all higher $2^{N}$-ions, $N \geq 4$, almost
instantly squelched the ambitions of number theorists -- with such
thoroughness, that the next (32-D) numbers were never even given a
name, much less investigated seriously. Yet these Pathions, as we'll
term them (as in ``pathological,'' which is what such entities were
deemed) are precisely those which hold the key to underwriting (and
algebraically generalizing) fractals with Number Theory.

There is much more than irony in the fact that fractals themselves
were similarly dismissed as ``monstrosities'' -- only to be later
tamed and rendered commonplace by Benoit Mandelbrot. [2] For, thanks
to their attribute of being scale-free (as statistical
distributions), fractals buttress the recent theory of
\textit{complex networks} [3] -- which theory can be given its own
Number Theoretic basis (making it henceforth less reliant on
empirical methods) by dint of ZD ensembles. Such a network-enabling
meta-fractal (each of whose infinite points is associated with one
among an infinite number of orthogonal lines of ZDs) we call a Sky
-- for it is where Box-Kites fly. These latter ensembles of 6
orthogonal planes, each represented by a vertex on an octahedral
frame, are the basic molecules of ZD structures. Their workings (and
manner of embedding in Skies, which we'll visualize via
spreadsheet-like ``emanation tables'' or ETs) will be our first
concern.  Once we understand their interrelations, we will be able
to generate the simple rules that use the bit-string of a Box-Kite
ensemble's integer signature to fill or hide ET cells, in a manner
implying approach to the ``fractal limit.''

Our argument, like Caesar's Gaul, naturally divides into three
parts.  In this Part I, our focus will be on the octahedral Box-Kite
representation, which can organize all the seminal concepts of
Zero-Divisor Theory in the 16-D context in which it first emerged,
and end with a first foray into 32-D, fabricating the small set of
fundamental tools we'll need to keep expanding our agenda
recursively into 64, 128, 256, then infinite dimensions.  In Part
II, we introduce the spreadsheet-like means of representing
ensembles of Box-Kites -- the Emanation Tables alluded to last
paragraph -- which lead to a graphics more suggestive of colored
quilts than any regular figures Euclid might have known.  Their
resemblance to classical fractals will soon be evident, and
culminate in our construction of the simplest (``Whorfian'') Sky,
whose basic building blocks only first emerge in 32-D. Finally, in
Part III, we will construct the component parts of the algorithm
which lets us build an infinite-dimensional meta-fractal Sky that is
unique to any integer $> 8$ and not a power of 2 (which is how we
can justify calling our study a branch of Number Theory -- but one
based on bit-string-suggestive placeholder structures, with powers
of 2 playing a role roughly akin to that of prime numbers in the
quantity-fixated contemplations of which the Greeks were so fond).
Here, we conclude with what literally boils down to``cookbook
instructions'' for building arbitrary Skies, allowing for NKS-style
movements between them -- the Recipe Theory suggested in our
abstract.

\section{Preliminaries:  CDP; What Zero-Division Means}
We use the most general CDP labeling scheme:  units are subscripted
from 0 (Reals) to $2^{N}-1$; the product (up to sign) of any two has
subscript equal to the XOR of the producing terms' subscripts.  Call
{\bf G} the index of the CDP generator of $2^{N+1}$-ions from
$2^{N}$-ions.  Then ${\bf G} = 2^{N}$ ( $= 8$ for Sedenions), and we
have Rule 1:  for $2^{N}$-ion index $L < {\bf G}$, $i_{L}\cdot i_{G}
= +i_{(G + L)} \equiv +i_{(G ~ xor ~ L)}$. Once we have at least one
associative triplet (trip) to work with, we can invoke Rule 2:
writing in cyclically positive order (CPO) so positive units
multiplied left to right always result in another, creating a new
trip by adding ${\bf G}$ to \emph{2} terms of another
\textit{reverses the order}: ${\bf(a, b, c)} \rightarrow {\bf(a, c +
G , b + G)}$. These 2 rules (plus a Rule 0 which says trips in
$2^{N}$-ions remain trips in all higher $2^{N + k}$-ions, $k > 0$)
completely describe trip-making, hence CDP.

\smallskip
\textit{Examples.}  Appending a unit $i_{2}$ to the standard Complex
plane, Rule 1 gives us $i_{1} \cdot i_{2} = +i_{3}$ -- the usual
$i,j,k$ of Quaternions ($2^{2}$-ions).  Performing one further
iteration with ${\bf G} = 4$, Rule 1 applied to the singleton Q-trip
(1,2,3) yields ${\bf G - 1} = 3$ O-trips $(1,4,5); (2,4,6);
(3,4,7)$. By Rule 2, we also derive 3 more:  $(1,\textbf{2,3})
\rightarrow (1,\textbf{7,6})$; $(2, \textbf{3,1}) \rightarrow (2,
\textbf{5,7})$; and, $(3, \textbf{1,2}) \rightarrow (3,
\textbf{6,5})$.  Altogether, then, CDP says Octonions contain $1 + 3
+ 3 = 7$ O-trips -- the \textit{trip count} for $N=3$, or just
\textit{$Trip_{3}$}, a number we derive independently by simple
combinatorics: to form an associative triplet, pick one unit from
those available, then pick a second from those remaining, and then
divide by 6 to allow for all the permutations that could have led to
the designated triplet being fixed by 2 selections. (Hence, the
Sedenions, with $N = 4$, have $(2^{N}-1)(2^{N}-2)/3! = 35$ S-trips;
Pathions, with $N = 5$, have $155$ P-trips; $Trip_{6}$ = 651, and so
on.)

\smallskip
\textit{Remarks.} The unique Q-trip can be rotated to bring its
units' indices into ACO (ascending counting order) as well as CPO:
 $(1,2,3)$.  2 of the 7 O-trips, though --  $(1,7,6)$ and $(3,6,5)$ --
 are ``bad trips.''  Inspection of our examples indicates that of the 3
 new trips resulting from applying Rule 2 to the Q-trip, only
 one -- O-trip $(2,5,7)$ -- will be ``good'' in turn.  Corollarily,
 starting with the Sedenions, Rule 2 will produce 2 ``good'' and 1
 ``bad'' trip for each ``bad trip'' fed into it: for instance, $(1,7,6) \rightarrow
 (1, 14, 15); (7, 9, 14);(6, 15, 9)$.  Simple algebra readily
 persuades us that $Trip_{N+1} = 4 \cdot Trip_{N} + 2^{N} - 1$; further, that ``bad trips'' among the $2^{N+1}$-ions = $2 \cdot
 Trip_{N}$, and ``good,'' this number augmented by the Rule
 1  contribution, or $2 \cdot Trip_{N} + {\bf G} - 1$.  Hence, as $N \rightarrow
 \infty$, good and bad trips approach parity.

\section{Getting Started:  Box-Kite Exclusion Rules}

Treat the 8 triangles of an octahedron as like a checkerboard, with
those of same color meeting only at vertices.  But instead of
colors, envisage 4 are Sails (jibs X, Y, Z, W, made of mylar, say),
while those on the faces opposite them are Vents (through which the
wind blows). The 3 pairs of opposite vertices are separated by
\textit{struts} -- wood or plastic doweling in real-world kites.
Such \textit{strut-opposites} are the \textit{only} pairings of
vertices which are \textit{not} mutual ``divisors making zero''
(DMZs): traversing any of the 12 edges always ``makes zero,'' as the
ends of each host two (mutually exclusive) DMZ pairs.

\begin{figure}
\centering
\includegraphics[width=.6\textwidth] {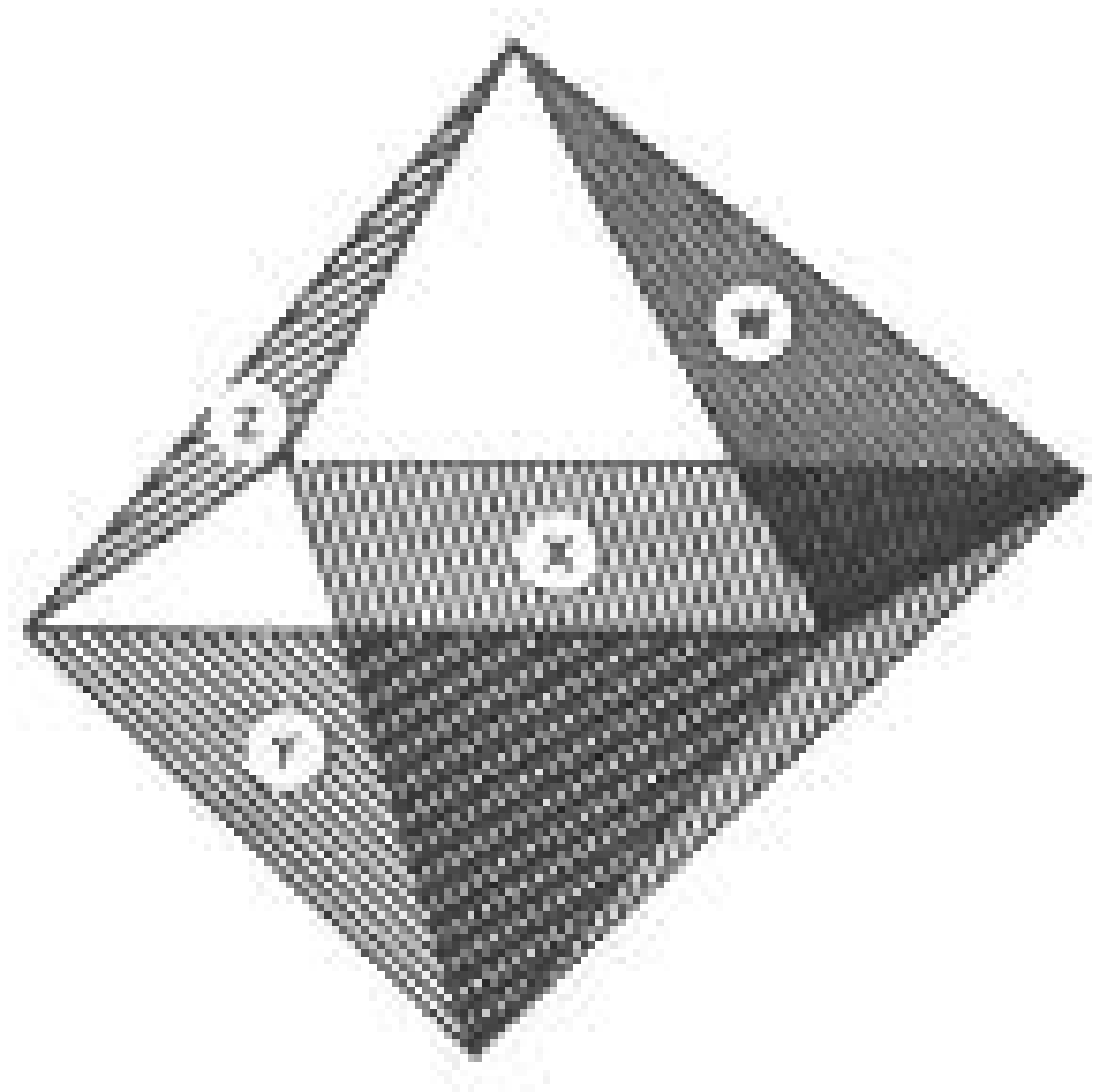}
\end{figure}

For this to make sense, we must explain what we mean by ``traverse''
in this context. At each vertex are two imaginary units of different
index (one, as we'll show, always with index less than {\bf G} -- A
or B, say -- while the other is always of index greater -- H or J).
The edge between any two Box-Kite vertices which aren't strut
opposites \textit{always} connects two such dyads whose possible
products include zeros.  Their general product consists of four
components, created by simple term-by-term multiplication. With
proper sign choices (a crucial matter to which we'll return
shortly), the product of the units indexed B and H will exactly
cancel that of the units indexed A and J; likewise, the product of
those indexed A and B will exactly cancel that of units indexed H
and J. (Note, in this last case, the high bits of H and J cancel by
XOR, whereas A and B, being less than {\bf G}, lack these high bits
to begin with.)  Also, we say \textit{two} sets of DMZ pairs framing
each edge, since each vertex represents one among 6 orthogonal
\textit{planes}, each spanned by a pair of hypercomplex axes, whose
\textit{diagonals} house all (and only) ``primitive'' ZDs.

For the case just given, we write $ (i_{A} + i_{H})$ and $( i_{A} -
i_{H})$, with one of these providing 2 of the 4 terms in the zero
product made with a $(B, J)$-indexed diagonal.  We will see in
Theorem 3 that 2 of the 4 possible choices of diagonals to multiply
together will always have zero product, and determining which
pairings ``make zero'' is where signing comes in.  Put another way,
since CDP-generated imaginaries square to $(-1)$ tautologically, a
ZD is ``primitive'' only if comprised of the sum or difference of
\textit{two} such imaginaries. (DMZ also signifies ``\textit{dyads}
making zero.'') ``Traversing'' an edge, then, entails multiplying
arbitrary points on suitably matched diagonals, one chosen from each
terminating dyad, allowing some choice (since the order in which we
multiply them has no effect on whether we get a zero or not) as to
which end we think of as the starting point.

The reason for this image will become clear when we consider
chaining such traversals: certain sequences of vertex-to-vertex
multiplication effectively take us on closed tours of some subset of
the Box-Kite's pathways, in such a way that we always ``make zero''
at each step in our circuiting. There are only a few such closed
circuits, falling into a small number of types. These have
surprisingly rich and distinctive properties as we shall see.  Our
focus here, though, will be primarily on the 3-vertex, 6-cyclic,
kind, which we've already given a name:  the Sail.  (The third
vertex is \textit{emanated} by the first two ZDs, since the zero the
latter make is in fact the sum of two oppositely signed copies of
this third dyad's terms.  What makes the Sail an algebraic whole is
the fact, shown in Theorem 5, that \textit{any} of its 3 edges, when
traversed, emanates the third vertex in this sense.)

As the general reader is likely not a trained algebraist, some
connection between such notions and those more familiar to
computer-based thinking are appropriate.  We note, then, that no
notion is more fundamental than that of a solution, and the
algebraic study of solution spaces can be seen as fundamentally
devolving upon zero divisors and imaginary quantities
simultaneously.  The very origins of the algebraic theory of
equations demanded imaginaries as intermediary terms, as Cardano
found when he discovered the general form of the cubic equation in
the Renaissance.  And as is well known, the intractibility of
quintic equations to extensions of such methods led to the classic
group-theory-spawning studies of Abel and Galois. And here, we find
the symmetry groups of equations devolving upon the structure of the
zeros they contain.

Taking the step from calling the two terms $(x - b)(x - c) = 0$ in a
quadratic equation, say, mere factors, to legitimate zero-divisors
in our sense here, is one in keeping with historically recent
paradigm shifts in physics.  The move from point particles to 1-D
strings and multi-dimensional n-branes is roughly parallel to our
own shift from zero-divisor \textit{points} to diagonal
\textit{lines} and exotica like the Seinfeld hyperplanes of our
November, 2000, paper. And just as one gets from points to lines by
stringing an indefinite number of the former together in an
ensemble, we will expand the singleton imaginary fundamental to
Galois' thinking to the infinite array of such units afforded by
CDP.

As CDP will be our starting point when our argument resumes, we
allow ourselves one final comment before returning to it.  Galoisian
zeros offer a basis for modern representation theory, providing the
starter kit (once embedded in Sophus Lie's continuous-groups
treatment) for handling all linear forms (and even the
``elementary'' singularities of nonlinear kind comprising Ren\'e
Thom's now-classical Catastrophe Theory). But the representation of
the nonlinear, we claim, requires methods in the general case like
those we adumbrate (and, to some small degree, actually implement)
herein.  And, knowing nothing more than the CDP rules, we can
exclude some very general cases:

\medskip \noindent {\small Theorem 1.}  Neither dyad can contain only
indices $< {\bf G}$.

\smallskip
\noindent \textit{Proof.}  If both dyads contain only Octonions,
their product will contain only Octonions, which have no ZDs; by
induction, no ZDs new to the $2^{N}$-ions generated by {\bf G} can
emerge from products of terms with indices $< {\bf G}$.  If only one
dyad has both terms' indices $< {\bf G}$, the term $\geq {\bf G}$ in
the other will have unequal products with the first dyad's two
terms. But, by XOR, these will also be the only terms with indices
$\geq {\bf G}$, preventing cancelation of either. {  }$\blacksquare$

\medskip
\noindent {\small Theorem 2.}  For any $2^{N}$-ions, the unit
indexed by {\bf G} cannot belong to a ZD dyad.

\smallskip
\noindent \textit{Proof.}  Given trip $(a,b,c)$ and binary variable
$(sg) = \pm 1$, write $(i_{a} + i_{G})$ below, and $(i_{b} + sg\cdot
i_{G + c})$ above, and multiply term by term.  The products by
$i_{G}$ are $-i_{G + b}$ and $sg \cdot +i_{c}$ (Rule 1). Products by
$i_{a}$, however, are $+i_{c}$ (given) and $sg \cdot +i_{G + b}$
(Rule 2). Depending on $sg$, one or the other, but not both, pairs
will cancel.  {  }$\blacksquare$

\smallskip
\noindent \textit{Remark.}  When $N$ increases, former generators
form ZDs. Starting with the Pathions, though, a condition we term
\textit{carrybit overflow} can arise, wherein indices sufficiently
large to play different roles for different $N$ cause the breakdown
of zero-division behavior in many cases where we would expect it,
and it is precisely because of this that the meta-fractal Skies will
arise.

\medskip
\noindent {\small Theorem 3.}  Suppose two dyads each contain one
unit whose index $> {\bf G}$, written in upper case, and another
whose index $< {\bf G}$, written lower case.  Then if $(a + A)\cdot
(b + B) = 0$, $(a + A) \cdot (b - B)$ does not, and vice-versa.

\smallskip
\noindent \textit{Proof.}  Given the zero result, the aB and Ab
terms, however signed, must cancel, as must AB and ab. But changing
only \textit{one} sign in one or the other dyad reverses sign on one
unit only in each pair. {
 } $\blacksquare$

\smallskip
\noindent \textit{Remark.} Clearly, multiplying either dyad by a
real scalar $k$ has no effect on ZD status, which is why these ZDs
are best thought of as \textit{diagonal line elements}.  Represent
ZD diagonals $k(x + X)$ and $k(x - X)$ in the obvious graphical
manner as $\verb|(X,/)|$ and $\verb|(X,\)|$ respectively.  It is
easy to show these corollaries:

\begin{center}
$\verb|(A,/)| \cdot \verb|(B,/)| = 0 \Leftrightarrow
\verb|(A,\)| \cdot \verb|(B,\)| = 0$

$\verb|(A,/)| \cdot \verb|(B,\)| = 0 \Leftrightarrow \verb|(A,\)|
\cdot \verb|(B,/)| = 0$
\end{center}

\medskip \noindent {\small Theorem 4.}  Two ZD diagonals spanned by
the same imaginary units cannot be DMZs.

\smallskip
\noindent \textit{Proof.}  $(a + A)(a - A) = +Aa + 1$ on top, and
$-1 - aA = -1 + Aa$ on the bottom.  But $2Aa \neq 0$. {
 }$\blacksquare$

\smallskip
\noindent \textit{Remark.}  This is the exact opposite of the
(non-CDP) ZDs found in quantum mechanics (QM), where the mutual
annihilation of the idempotent projection operators $\frac{1}{2}(1
\pm m)$, $m$ a Pauli spin matrix (so that $m^{2} = +1, ~ m \neq \pm
1$), serves to define observables. There is a sense, however, in
which the QM case is a degenerate form of the behaviors being
examined here: this would take us off topic, but those interested in
such questions should see [4].

A final exclusion rule, alluded to above:  ZDs living in planes
associated with strut-opposite vertices can never be DMZs.  This is
a deeper result than those just presented, and is a side-effect of
the general production rules, to which we turn next.

\section{Completing the Schematic:  Box-Kite Production Rules}
In the Sedenions, the most general production rule allows anything
not excluded by the rules above:  pick any Octonion of index $O$ (7
ways), then any pure Sedenion with index $S > {\bf G}$ not the XOR
of $O$ with {\bf G}; the resulting dyad spans one of 42 planes --
the 42 Assessors of [5], where the other production rules are also
described for the first time.  For higher $2^{N}$-ions, restrictions
due to carrybit overflow arise, but the approach just limned still
generates all \textit{candidates} for ``primitive ZD'' status.

From this rule to Box-Kites is but a short distance:  for, since
each Sedenion Assessor will have an Octonion in its dyad of units,
clearly we can arrange their 42 planes in 7 clusters of 6, each of
which excludes a different Octonion from its set.  Call the Octonion
index in each dyad (and, more generally, the $2^{N}$-ion index $<
{\textbf G}$) the L-index, and its partner the U-index, and write
their respective units $i_{L}$ and $i_{U}$.  Since each Octonion is
involved in precisely 3 O-trips, the 6 L-indices in each cluster can
be uniquely fixed in 2 steps:  first, assign each pair whose XOR is
the excluded index to one of the 3 struts, which is necessary and
sufficient to guarantee that their product does not appear in the
set.  (From what we've already seen, even without specifying what
their U-indices are, this is tantamount to guaranteeing that no
pairing of their Assessors' diagonals can be DMZs.)  Call the
excluded index {\bf S}, for ``Strut Constant'' and ``signature.''
The last parenthetical remark tells us it somehow fixes the
U-indices of the Assessors in its cluster.  We will see how in the
next section.  For now, the more immediate problem of fixing
L-indices leads us to perform our second step.

Orient the strut-opposite L-indices so that {\bf S} is the central
label in the PSL(2,7) triangle which represents the smallest finite
projective group, also called the Fano plane.  Its 3 vertices, 3
sides' midpoints, and center each host an Octonion index, in 7
oriented lines of 3 points each, one per O-trip (in 480 possible
configurations). Then, assign the L-index trip in the circle
threading the sides' midpoints to vertices $A$, $B$, $C$ on a
Box-Kite, in CPO order. The strut-opposite L-indices will then
appear as the vertices which terminate the lines extending from
their L-index partners through {\bf S}.  We label them, in ``nested
parentheses'' order, $F$, $E$, $D$ respectively.

By removing the central label, four of the 7 lines in PSL(2,7)
remain connected:  that of the central circle just discussed, whose
L-indices form the triplet $(a,b,c)$; and, the 3 along the
triangle's sides, which the diagram tells us have CPO labels
$(a,d,e)$, $(f,d,b)$, and $(f,c,e)$.  We will show next section how
these form the bases of the two types of Sails -- the singleton
Zigzag at vertices $(A, B, C)$, and the three Trefoils which each
share a different vertex with the Zigzag -- and demonstrate their
special properties, which hold for all Box-Kites in all
$2^{N}$-ions. For now, we introduce the other production rules, the
first of which motivates the very idea of a Sail.

\pagebreak

\noindent {\small Theorem 5.}  If Assessors $(U,u)$ and $(V,v)$ form
DMZs, then a third Assessor $(W,w)$, forming DMZs with both, can be
created by setting the index of $W$ to $|U ~xor~ v| = |u ~xor~ V|$,
and that of $w$ to $|u ~xor~ v| = |U ~xor~ V|$.

\smallskip
\noindent \textit{Proof.}  Without loss of generality, assume
L-index trip $(u,v,w)$ is CPO, and the dyads of the given DMZ have
same internal signing: that is, $(u + U) \cdot (v + V) = 0$. We then
have $u \cdot v = - U \cdot V$; but $(u,v,w) \rightarrow v \cdot w =
+u$, so $w = - U \cdot V$ as well.  Using $sg$ as before, write $W =
sg \cdot [ + Uv ] = sg \cdot [ - uV ]$. Consider $(v + V) \cdot (w +
sg \cdot W)$, the former written under the latter. (The same basic
argument applies if we choose $(u + U)$ instead and write it on
top.) Multiplying left to right by $v$, we have $+u$ and $sg \cdot
(v \cdot Uv) = sg \cdot ( -v \cdot vU ) = sg \cdot ( +U )$. Next, $V
\cdot w = V \cdot ( - UV ) = V \cdot ( +VU ) = - U$. And $V \cdot sg
\cdot W = V \cdot sg \cdot ( -uV ) = sg \cdot ( - u )$. Set $sg =
(-1)$, and all terms cancel. { } $\blacksquare$

\smallskip
\noindent \textit{Remark.}  Two Assessors in the same Box-Kite
\textit{emanate} a third whose L-index completes the L-trip implied
by the other two. (The zero they make is also a sum of oppositely
signed copies of the emanated Assessor.)  Such a triad defines a
\textit{Sail}.

\smallskip
\noindent \textit{Caveat. } The proof as given is complete within
the Sedenions, which display alternativity:  $(aa)\cdot b = a \cdot
(ab)$. Generalizing to higher $2^{N}$-ions, however, requires a
further result we'll derive later, namely: all three index-triples
involving one lowercase and two uppercase letters, each from a
different Assessor in a Sail, are trips, so that even
\textit{associativity} among units combined as required in our proof
is guaranteed.

\medskip
\noindent {\small Theorem 6.}  If $(U + u)$ and $(V + sg \cdot v)$
are DMZs, then so are their \textit{twist products} $(U + sg \cdot
v)$ and $(V - u)$.

\smallskip
\noindent \textit{Proof.}  By direct multiplication, the DMZ given
requires $U \cdot v = - u \cdot V$ and $u \cdot v = - U \cdot V$.
The proposed DMZ resulting from the twist would require $V \cdot v =
u \cdot U$, and $u \cdot v = V \cdot U$.  But the second requirement
is identical to what we know is true in the given DMZ.  As for the
first requirement, multiply the terms of the other known relation on
the right by $u \cdot v = w$ to obtain it. { }$\blacksquare$

\smallskip
\noindent \textit{Caveat added in V3. } The proof as given is
complete within the Sedenions; however, it doesn't hold universally
for the Pathions and beyond, when strut constants can exceed 8 -- a
side-effect, and necessary correlate, of the emergence of the
``meta-fractal'' structures that will concern us in the second and
third installments of this argument.  Perhaps surprisingly, this
proves to be a virtue: the lanyards based on circuiting squares, not
triangles, are the natural vehicles with which to see twist products
deployed, and it is thanks to the twists which do \textit{not} yield
ZDs that we can study transformations between legitimate box-kites
and the ``broken'' kind lacking ZD pathing which comprise
meta-fractals. Such topics, however, will not be broached in any of
the three parts of this monograph; rather, they require their
background to be dealt with, and hence will provide the subject
matter for their immediate sequel.

\smallskip
\noindent \textit{Remark.}  A few months after [5], and quite
independently of it, Raoul Cawagas obtained the same listing of
primitive ZDs in the Sedenions.  He obtained it by programmatic
exploration, using his Finitas software, of \textit{loops} (which do
for nonassociative algebras what groups do for their more orderly
cousins). [6]  In our correspondence, it became clear that twist
products in fact connect ZDs which, while necessarily in different
Box-Kites, are always in the same loop. Each of a Box-Kite's 4 Sails
``twist'' to different loops, while each Cawagas loop meanwhile
partitions into 4 sets that correspond to Sails in different
Box-Kites, indicating an interesting duality between the two quite
differently framed approaches and their objects.  For the surprising
symmetries inherent in twist product patterns, see the discussions
surrounding the Royal Hunt and Twisted Sister diagrams in [7, pp.
14-16].  See [5, pp. 11-21] for full exposition of all the
production rules being discussed here, and pp. 15-16 for the twist
product rule in particular.

\smallskip
\noindent \textit{Coming Attractions.}  Once the workings of the
Strut Constant are understood, the reader will find it trivial to
prove that opposite sides of any of the 3 orthogonal squares
comprising the Box-Kite's octahedral framework twist to the same
Box-Kite, while adjacent sides twist to different ones.  (The trick
is to realize that, for all Assessors in any given Box-Kite, the XOR
of their dyads' indices is ${\bf G + S}$).  This suggests that the 3
squares, which we will call Catamarans (Tray-Racks in earlier
papers), are as fundamental as the 4 Sails, but for the different
purpose of long-distance navigating; we conventionally name them for
the tall masts planted orthogonally on the transverse frames joining
their hulls:  $AF$ for Catamaran $BCED$, for example. Catamarans
also come in pairs:  as will also be clear when edge-signs are
described, tracing the perimeter of a Catamaran only engages 4 of
the 8 ZD diagonals in the Assessor planes being passed through,
whereas tracing a Sail twice will engage all 6 residing in its 3
Assessors. Catamarans, then, with their square perimeters and coming
in pairs, could easily be represented by traditional Box-Kites, as
opposed to our surreal octahedral ones.

\medskip
\noindent {\small Lemma.}  As every L-index appears in two distinct
L-trips in a given Box-Kite, likewise every Assessor appears in two
distinct Sails.

\smallskip
\noindent \textit{Remark.}  Full XOR-based calculation requires
knowledge of the U-indices, which requires, in turn, a better
understanding of the special role played by the Strut Constant --
which will be our focus in the next section.

\section{ Strut Constants, Sails, and ``Slipcover Proofs'' }
The debt to the groundbreaking work of R. Guillermo Moreno is
evident:  his 1997 preprint (appearing in a Mexican journal the
following year [8]) revivified the long-neglected study of ZDs.
Moreno's now-classic paper culminated in a formal demonstration that
the algebra of ZDs in the Sedenions was homomorphic to the Lie
algebra $G_{2}$, itself the automorphism group of $E_{8}$, in turn
represented by the Octonions -- a result seemingly of interest only
to physicists. Yet he also pointed out that this same $G_{2}$ is the
derivation algebra of $2^{N+1}$-ions from $2^{N}$-ions, implying a
deepset recursiveness in ZD structuring.  We underscore that the
Fano Plane diagram for Octonion labeling is likewise a potentially
recursive tool:  the 15 points and 35 lines, with 7 lines through
every point and 3 planes in each line, of the 3-D ``projective
tetrahedron'' that schematizes Sedenion labeling, is further cobbled
together from 15 planes isomorphic to the $PSL(2,7)$ triangle.
Higher labeling structures of $2^{N}-1$ points, $Trip_{N}$ lines,
and so on, are readily imagined (albeit impossible to visualize).
Our next results will be derived with virtually no tools beyond CDP
and PSL(2,7).  We begin with a concretization of the discussion
which introduced the latter last section.

Place the indices of the Octonions so as to reflect their
construction via CDP from Quaternions:  in the center of an
equilateral triangle, place ${\bf G} = 4$; with apex extending above
12 o'clock, place $(1, 2, 3)$ at the midpoints of the left, right,
and bottom sides, oriented clockwise, at 10, 2, and 6 o'clock in
that order. Now its projective line, the only one drawn as a circle,
represents the Rule 0 trip; lines oriented from the midpoints to and
through the center are the Rule 1 trips; the sides, meanwhile,
constitute the Rule 2 trips, and are all oriented in a clockwise
manner, paralleling the flow along the Rule 0 circle.

Now replace $(1, 2, 3)$ with the letters $(u, v, w)$, and replace
the indices at the angles with the symbols $(u_{opp}, v_{opp},
w_{opp})$, as these are indeed $u$, $v$ and $w$'s strut-opposites
with respect to the central index, which we now replace with the
symbol {\bf S} for Strut Constant.  The contents of the labels are
indeterminate, but their choices must conform to the network of
flows:  clockwise along edges and central circle, and pointing
inward along the 3 lines bifurcating sides, then angles.  In our
initial setup, the Octonion generator is indeed the {\bf S} for the
Sedenion Box-Kite implicated.  In a manner suggestive of how one
pulls on and removes slipcovers from upholstery, we claim that this
flow pattern is not only uniquely linked to representations with the
central circle comprising the L-indices of the Zigzag Sail; but that
one can, without inducing any tearing in the network of flows,
``pull'' any Octonion index into the center, thereby representing a
Zigzag-centric depiction of the Box-Kite whose ${\bf S}$ is that
index.

PSL(2,7) has the (for our purposes) highly useful property that, in
order to obtain a correct multiplication scheme for the 7 units, it
suffices to place $(u, v, w)$ on any line, mark any of the remaining
4 points ${\bf S}$, and then pass through it to designate the last
trio of points as opposites of the first.  Swapping in one point for
another, then, is easy.  The symmetries are so complete that we can
take any point we please and make the case for all, so let's take 7
in the ${\bf S} = 4$ diagram.  Rotate the line containing 4 and 7
once in the downward direction to keep orientation:  the vertical
line, in top-down order, now reads $(3, 7, 4)$ instead of $(7, 4,
3)$. The central circle, if its clockwise orientation is to be kept,
must be CPO starting with the left midpoint and ending in the 4 now
on the bottom.  There are just 2 possibilities -- $(6, 2, 4)$ and
$(5, 1, 4)$ -- but the former would reverse the flow on the right
flank, reading $(1, 2, 3)$ from the lower right to the apex. As is
readily checked, $(5, 1, 4)$ does exactly what we want -- and, as
our table of Sedenion Assessors and Strut Constants makes clear,
that is in fact the Zigzag L-trip for the ${\bf S} = 7$ Box-Kite.
The 3-fold symmetry of the flows lets us choose any of the central
circle's 3 nodes to house the $1$ (which is the L-index of the A
Assessor in the Box-Kite in question) without effecting the pattern.

\begin{figure}
\centering
\includegraphics[width=1.0\textwidth] {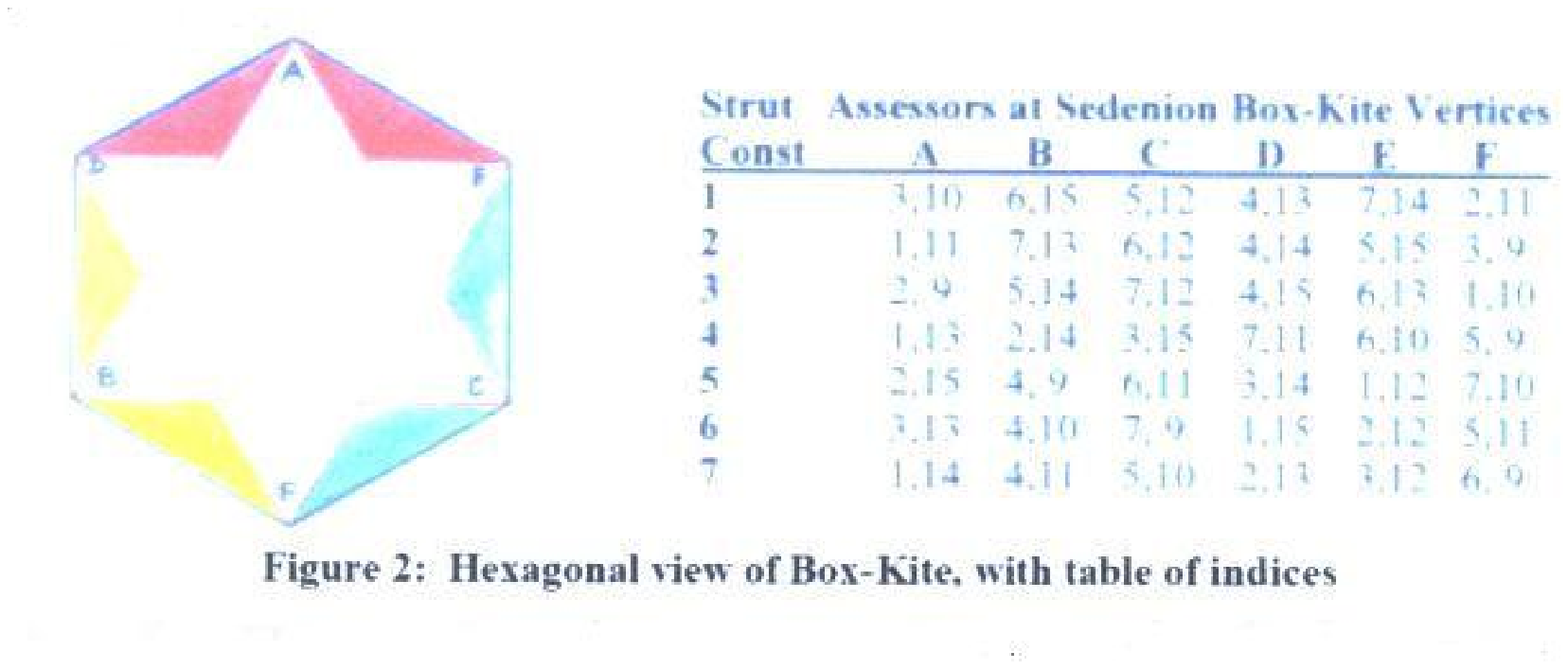}
\end{figure}

But suppose we wish to swap two \textit{lines} instead of nodes: one
of the Trefoil Sails on the sides for the central Zigzag, say.  This
can be done, but at the cost of changing the flow structure. Keeping
the convention of mapping $(a,b,c)$ to $(u,v,w)$, the Trefoils
$(a,d,e)$, $(f,d,b)$ and $(f,c,e)$ now get rewritten as $(u,
w_{opp}, v_{opp})$, $(u_{opp}, w_{opp}, v)$ and $(u_{opp}, w,
v_{opp})$.  If we insist on maintaining the clockwise orientation of
the central circle after replacing its labels with one of the 3 sets
just written, the flows on the 2 sides not corresponding to the
original $(u,v,w)$ circle will be counterclockwise, and the 2 rays
through the center starting from midpoints suffixed ``opp'' will now
lead away from, not into, the angles.  ``Slipcover tugs'' on the
fabric of \textit{this} network will send a Trefoil linked to one
${\bf S}$ value into one linked to a different ${\bf S}$ -- but the
lack of the Zigzag L-trip's threefold symmetry means it typically
will be a \textit{different} Trefoil:  if one tugs $2$, say, into
the center, the $(a,d,e)$ Trefoil of the ${\bf S} = 4$ Box-Kite is
transformed into an $(f,d,b)$.  The side-effect of Trefoils' broken
symmetry explains why, in the table provided above, Zigzags are so
regular (each O-trip appears as one \textit{exactly once}), while
Trefoils behave so erratically (each O-trip appears thrice as a
Trefoil, but \textit{never} in an evenly distributed manner; in
fact, $(3,6,5)$ makes all its Trefoil appearances as $(f,c,e)$,
while $(1,2,3)$ \textit{qua} Trefoil shows up only as $(a,d,e)$).

We mentioned above that PSL(2,7) diagrams can be used to map not
merely the interconnections of Octonion labelings, but -- given the
recursive modularity of all higher $2^{N}$-ions -- heptads among the
Sedenions and beyond.  We have used, for instance, ``stereo Fano''
diagramming to navigate and notate the structures found in the 32-D
Pathions.  [9]  For the present, though, we'll need nothing so
exotic:  we'll just want to drop ${\bf G}$ into the center of
diagrams already seen, not to replace ${\bf S}$, however, but to
supplement it.  The reason for this is implicit in this critical
theorem.

\medskip
\noindent {\small Theorem 7.} The U-index of any Assessor is the XOR
of its L-index with ${\bf G + S}$.  For Zigzag Sails, all 3 edges
connecting Assessors conform, in the notation of the Remark to
Theorem 3, to this pattern:

\begin{center}
$\verb|(M,/)| \cdot \verb|(N,\)| = 0 \Leftrightarrow \verb|(M,\)|
\cdot \verb|(N,/)| = 0$
\end{center}

\noindent For Trefoil Sails, the 2 edges including the Assessor
shared with the Zigzag are governed, instead, by the other pattern
from the same Remark, namely:

\begin{center}
$\verb|(M,/)| \cdot \verb|(N,/)| = 0 \Leftrightarrow \verb|(M,\)|
\cdot \verb|(N,\)| = 0$
\end{center}

\smallskip
\noindent \textit{Proof.}  Place $(u,v,w)$ on the central circle of
a PSL(2,7) diagram, with all 7 lines in zigzag flow orientation.
Place ${\bf G + S}$ in the diagram's central node, calculate the
angle nodes, and readjust flows by Rule 2:  the trips on the lines
through the center are now oriented to point from angles to
midpoints, and are written thus: ${\bf (G + u_{opp}, G + S, u)}$;
${\bf (G + v_{opp}, G + S, v)}$; ${\bf (G + w_{opp}, G + S, w)}$.
The side trips now flow \textit{counter}clockwise.  In CPO, from the
top, we have these: ${\bf (G + w_{opp}, u, G + v_{opp})}$; ${\bf (G
+ v_{opp}, w, G + u_{opp})}$; ${\bf (G + u_{opp}, v, G + w_{opp})}$.
We assert $U = {\bf (G + u_{opp})}$ and $V = {\bf (G + v_{opp})}$.
Multiply $(u + U) \cdot (v - V)$, the indices assumed to be
subscripts of implied imaginary units.

\begin{center}
$+ v - (G + v_{opp})$

\underline{$+ u + (G + u_{opp})$}
\smallskip

$+ (G + w_{opp}) {  } - w$

\underline{$ + w {    } - (G + w_{opp})$}

\smallskip
$0$
\end{center}

\smallskip
Now place $(u, w_{opp},v_{opp})$ on the central circle of a
trefoil-flow-oriented PSL(2,7); plunk ${\bf (G + S)}$ down in the
center of it; multiply dyads as above, but with \textit{same} inner
sign for the two products including Assessor U, to complete the
proof. { }$\blacksquare$

\medskip
\noindent \textit{Corollary.}  As a side-effect of Theorem 7, we see
that indeed each Sail is comprised of 4 trips, thereby addressing
the Caveat to Theorem 5.  The Zigzag, for instance, has one L-trip,
$(a,b,c)$, and 3 U-trips where high bits are lost by XOR, written
$(a,B,C)$; $(A,b,C)$; $(A,B,c)$. Each, supplemented by the Reals,
makes a legitimate copy of the Quaternions.  Moreover, any Box-Kite
in any $2^{N}$-ion domain will have Sails with this property, it
being a feature of our Rule 2 and PSL(2,7)'s recursivity, rather
than anything specific to (hence subject to the limits of) Sedenion
space.  We merely use whatever ${\bf G}$ is appropriate in our
``slipcover proofs'' and choose a Rule 0 triplet to place in the
central circle; or, equivalently, pick whatever trio of
strut-opposites you care to use which satisfy the given ${\bf S}$
value. (There will be many to choose from as N grows large, so that
we speak more of Box-Kite \textit{ensembles} than Box-Kites as
individuals).

\section{Intermezzo:  Lanyards, Semiotic Squares, Emanation Tables, Sand Mandalas \dots}

Moreno's ${\bf G_{2}}$, as root system, has order 12 (14, if we're
physicists and add in the 2-D Cartan subalgebra).  The  Sedenions' 7
Box-Kites x 12 edges x 2 oriented flows along each $\rightarrow$
168, order of the second simple group, and 12 \textit{times} 14. The
integer 168 is also a ``magic number'' in Boolean function theory,
where it arises in the context of the Dedekind Problem as the number
of 4-variable monotone functions, with the same number of
complements. The Dedekind \textit{Problem} concerns determining this
rapidly growing number for arbitrarily large counts of variables.
The relevance of such abstruse topics to rule-based approaches like
cellular automata theory and Stephen Wolfram's \textit{New Kind of
Science} has been the focus of ongoing work by Rodrigo Obando [10],
[9, pp. 18-19], who has employed their toolkit to create a kind of
Periodic Chart of such rule-driven systems, highlighting which are
not just Class 1 (boringly homogenous outcomes), Class 2 (evolving
into simple periodic structures), or Class 3 (chaotic, aperiodic),
but truly complex Class 4 patterns, like Wolfram's famous Rule 30.
All of which has provided motivation to search ZD space for
connections to similar modes of complexity, driven by similarly
simple rules.  We will conclude, as our abstract suggests, with just
such a ``recipe theory'' with which to study infinite-dimensional
meta-fractals. We will prepare for that journey by taking leave of
simple Box-Kites in the Sedenions, and go so far in this section as
to find the first truly anomalous ensembles, with the ``wrong
number'' of Box-Kites, in the 32-D Pathions.

We first want to indicate that even simple Box-Kites are much richer
in structure than we have seen thus far.  The two kinds of Sails are
instantiations of abstract objects we call \textit{lanyards} --
cords worn around the neck to hold knives or whistles, or strands of
leather or plastic used to string beads and other gewgaws into charm
bracelets.  By chaining together DMZs, a lanyard threads the ZD
diagonals of multiple Assessors,  attaining closure by returning to
its starting point having touched all other ZDs in its ambit but
once.  They can be given signatures even more compressed than the
notation developed in the Remark to Theorem 3:  the Zigzag is a
6-cycle lanyard, which ``makes zero'' with all its Sail's diagonals,
in a sequence cyclically equivalent to the glyph inspiring its name:
\verb| / \ / \ / \ |.  The Trefoil has its own 6-cycle signature:
\verb| / / / \ \ \ |.  (The Catamaran, meanwhile, has 4-cycle
signature \verb|/ / \ \|, which should clarify our earlier cryptic
comments concerning it.)

Joins between diagonals with opposite inner sign we designate either
with a [-] \textit{edge symbol}, else -- if availing ourselves of
color graphics -- draw in red.  Exactly half a Box-Kite's 12 edges
are red, consisting of the trio defining the Zigzag, and those
bordering the Vent on the face opposite it.  (See diagram on left of
Figure 2.) Those edges linking similarly inner-signed diagonals are
marked [+] or painted blue.  ``The Blues,'' in fact, are two other
6-cycle lanyards which contain all and only the blue edges,
corresponding to the perimeter of the hexagon in the prior image.
\textit{Two}, because one connects only \verb|\| diagonals, the
other only \verb|/|'s.  What might a simulation that switched
edge-signs (converting the Blues into Zigzags, perhaps?) correspond
to in real-world processes?

Another type of 6-cycle lanyard, the 2-Sail, 5-Assessor ``Bow-Tie''
whose doubly passed-through ``knot'' must be in the Zigzag, has deep
symmetry-breaking properties that may well play a fundamental role
in string theory [4].  In general, each distinct lanyard (and there
are others which the reader is invited to explore involving 4, 5,
and all 6 Assessors:  see Section III of [5] for more on some of
these) creates a distinctive context of dynamic possibilities,
rather than the semblance of an answer (as groups, too often, are
assumed to offer). And, unlike groups, and as one might expect in
$0$- rather than $1$- based algebraics, an identity element is
nowhere in sight, with strut opposites the closest thing to
inverses.  Rather, the zero resultant defining the transit of each
link in a lanyard's chain guarantees a certain design-hiding,
suggestive of the unfolding of the multitude of zeros packed in a
single seed-form or morphogenetic ``vanishing point.'' Put another
way, a theory of lanyards can be expected to provide some crucial
ideas in an area currently lacking many: the study of nonlinear
representations.

One area where a theory of nonlinear representations exists is in
semiotics, where the French school of Jean Petitot -- dually
inspired by Ren\'e Thom's Catastrophe Theory (CT) mathematics in the
1970's, and the structural rapprochement of semantics and syntax
essayed by Algirdas Greimas -- led to CT modeling of ``archetypal
nouns and verbs'' on the one hand, and a delving beneath Chomskian
notions of semantically impoverished recursiveness on the other.
Using the CT toolkit, Petitot has built models of Greimas' own
famous (among literary theorists, at least) Semiotic Square, [11] as
well as the equally celebrated (among, at least, the
structuralism-savvy with anthropological leanings) canonical law of
myths of Claude L\'evi-Strauss. [12] But both these models can be
embedded, in their turn, in ZD-based networks, which put the
recursiveness back where it belongs:  not at the level of meaning as
such, but rather in how meaning-rich nodes can be constellated in an
open-ended scale-free dynamic.

The key insights are two:  first, to see the entire program of
structural linguistics -- whose origins were contemporary with the
algebraic death-knell sounded by Hurwitz's Proof -- as a kind of CDP
\textit{manqu\'e}, thanks to Hjelmslev's indefinitely extensible
notion (seminal to Greimas' thinking) of form vs. content
\textit{double articulation}. [13]  The second insight is to find
semiotic substance in what would seem a \textit{lack} of substance
from our purview: by analogy with Shannon-Weaver self-correcting
error-code theory, where most of the bits are expended in the
error-trapping as opposed to the content being transmitted, we
assume the area of interest is precisely the ZD-challenged parts of
the ZD ensembles that convey the significances.  Specifically, we
find our meta-model of Petitot's CT reading of the Semiotic Square
in the 4 units whose invariant workings govern the ZD-free
strut-opposites:  the real unit (which never appears in primitive
ZDs), ${\bf G}$ and ${\bf S}$, and their XOR (which, given that
${\bf G}$ always has only one bit, always to the left of all of
${\bf S}$'s, is also their simple sum, henceforth to be called ${\bf
X}$).

The index set $\bf{\{0,G,S,X \}}$ defines the units of a Quaternion
copy, part of what Moreno terms the 4-D, ZD-free boundary of
Sedenion space, in which all Box-Kites, each with their own ${\bf
S}$ and ${\bf X}$, participate.  We argue at length in [7] that the
square with drawn diagonals -- $\boxtimes$ -- which Greimas made the
centerpiece of his thinking, can be modeled by placing ${\bf G}$ on
the diagonals, ${\bf X}$ on the horizontals, and ${\bf S}$ on the
verticals.  From the ZD perspective, the workings of the 3 modes of
binary relations (contradiction between a seme and its absence -- an
on/off bit -- along diagonals; contrareity of reciprocal
presupposition along the horizontals; complementarity of implication
along the verticals) can be formally encapsulated in the following
relations. Pick two strut-opposite Assessors, dubbing the Zigzag
dyad picked as $(Z,z)$, and the Vent dyad linked to it, $(V,v)$.
Then

\begin{center}
$v \cdot z = V \cdot Z = {\bf S}$

\smallskip
$Z \cdot v = V \cdot z = {\bf G}$

\smallskip
$V \cdot v = z \cdot Z = ({\bf G+S}) \equiv {\bf X}$.
\end{center}

These 3 relations were implicitly derived in passing while proving
Theorem 7.  Their importance (and not just for semiotic
applications) suggests we give them a special name:  the
\textit{Three Viziers}, to be shorthanded henceforth (and rewritten
in triplet form) as follows:

\begin{center}

${\bf VZ1:}\;\;(v,z,S); (V,Z,S)$

${\bf VZ2:}\;\;(V,z,G); (Z,v,G)$

${\bf VZ3:}\;\;(V,v,X); (z,Z,X)$
\end{center}

\smallskip
\noindent \textit{Caveat added in V3. } Only the second of these
three is truly universal; the first and third are likewise universal
\textit{only if unsigned}.  Beginning in the Pathions, one can find
Box-Kites where precisely two of the three struts have the relative
signing of vent and zigzag terms reversed, so that such struts have
their VZ1 L-index triplets, for instance, written ${(z,v,S)}$.
Implied in this is the existence of a second kind of box-kite, whose
nature and very existence only became apparent in work subsequent to
this monograph, once catamaran twisting was investigated.  (These
``Type II'' box-kites are in fact indistinguishable from the
``normal'' kind except in regards to such twisting, since struts do
not mutually zero-divide anyway, and hence their triplets'
orientational variations ``fly under the radar'' of the tools we've
built thus far.  They serve as ``middlemen,'' however, between the
``normal'' kind exclusively in evidence in the Sedenions, and the
ZD-free substructures which make up the body of the meta-fractals
we'll focus on next.) Skipping over such subtleties, fortunately,
has no effect on the \textit{outcome} of the general argument we
make in this monograph; a handful of steps in the \textit{proofs},
however, will require modification in Parts II and III.

\medskip
Getting from the atomic level of Greimas' Square, to the
small-worlds networking of literally thousands of mythic fragments,
all interconnected and mutually unfolded in the four fat volumes of
L\'evi-Strauss's magnum opus, would require some surgery on the
Square, to allow for a competition between contraries, say, to be
transmuted into an agreed-upon implication.  But the simplest such
transformative structures already require the Pathions, and hence
meta-fractal Skies to fly in.  Those are the vistas in front of ZD
theory, only some of which we will get to see herein. One such
surgery, corresponding to the type just suggested as an example,
will in fact provide the first recursive iteration leading to the
fractal limit. And, as we're about to see, it arises naturally,
virtually as a built-in feature, of the ZD framework.

As anyone seeing a ZD diagonal on an oscilloscope screen would
think, Assessors are ultimately patterns of synchrony and
anti-synchrony between their dyadic terms.  Assume Assessor-unit
synchrony breaks down, with dyads recast as strut-opposites.  Each
finds higher-index partners; the former ${\bf G}$ and ${\bf S}$ are
transformed into Pathion L-indices at ${\bf B}$ and ${\bf E}$, the
former ${\bf X}$ moves to center stage as the new ${\bf S}$.  14
Assessors (${\bf G} - 1$, minus 1 more to exclude the Strut
Constant) form into 3 Box-Kites sharing the ${\bf BE}$ strut.  It is
quite possible to draw this three-headed beast -- but the
many-headed hydras in higher $2^{N}$-ions will not prove so
amenable.  It is time to explore a different mode of representing
things, which means Emanation Tables.  And it is time to deal with a
new phenomenon whose crucial nature we've indicated more than once:
for values of ${\bf S} \leq 8$, the ensembles associated with each
such integer, contain seven, not three, Box-Kites, and are dubbed
\textit{Pl\'eiades} for this reason.   But the surprising graphical
sequence revealed as one thumbs through a \textit{flip-book}, whose
successive pages display \textit{emanation tables} linked to
progressively incremented ${\bf S}$ values from 9 through 15,
inspired a different name for these entities:

\begin{quote}
Viewed in sequence, these tables suggest the patterns made by
cellular automata; seen individually, they suggest nothing so much
as the short-shelf-life ``sand mandalas'' of Tibetan Buddhist
ritual, made by monks on large flat surfaces with colored sands or
powdered flowers, minerals, or even gemstones.  [14, p. 15]
\end{quote}

The second sweep of our argument [15] will begin with these two
themes as its ultimate agenda -- pursuing them in accordance with
our originally advertised intentions.  Once we have classified types
and characteristics of emanation tables through the 64-D Chingons,
we will have all the basic patterns needed for crafting (and
proving) the algorithmics of ``Recipe Theory'' in Part III. [16]

\section*{References}
\begin{description}
\item \verb|[1]|~ I. L. Kantor and A. S. Solodovnikov, \textit{Hypercomplex
Numbers:  An Elementary Introduction to Algebras} (Springer-Verlag,
New York, 1989)

\item \verb|[2]|~ Benoit Mandelbrot, \textit{The Fractal Geometry of Nature} (W.
H. Freeman and Company, San Francisco, 1983)

\item \verb|[3]| ~ Mark Newman, Albert-L\'aszl\'o Barab\'asi, and Duncan J.
Watts, editors, \textit{The Structure and Dynamics of Networks}
(Princeton University Press, Princeton and Oxford, 2006)

\item \verb|[4]| ~ Robert P. C. de Marrais, ``The Marriage of Nothing
and All:  Zero-Divisor Box-Kites in a `TOE' Sky'', in Proceedings of
the $26^{\textrm{th}}$ International Colloquium on Group Theoretical
Methods in Physics, The Graduate Center of the City University of
New York, June 26-30, 2006, forthcoming from Springer--Verlag.

\item \verb|[5]| ~ Robert P. C. de Marrais, ``The 42 Assessors and
the Box-Kites They Fly,'' arXiv:math.GM/0011260

\item \verb|[6]| ~ R. E. Cawagas, ``Loops Embedded in Generalized
Cayley Algebras of Dimension $2^{r}$, $r \geq 2$,''
\textit{International Journal of Mathematics and Mathematical
Science}, 28:3 (2001), 181-7.

\item \verb|[7]| ~ Robert P. C. de Marrais, ``Presto! Digitization'' arXiv:math.RA/0603281

\item \verb|[8]| ~ R. Guillermo Moreno, ``The Zero Divisors of the
Cayley-Dickson Algebras over the Real Numbers,'' \textit{Boletin
Sociedad Matematica Mexicana (3)}, 4, 1 (1998), 13-28; preprint
available online at arXiv:q-alg/9710013

\item \verb|[9]| ~ Robert P. C. de Marrais, ``The `Something From
Nothing' Insertion Point'',
http://www.wolframscience.com/conference/2004/presentations/materials/
\newline rdemarrais.pdf

\item \verb|[10]| Rodrigo Obando, ``Mapping the Cellular Automata
Rule Spaces,'' on http://
\newline wolframscience.com/conference/2006/presentations/obando.nb

\item \verb|[11]| Jean Petitot, \textit{Morphogenesis of Meaning},
(Peter Lang, New York and Bern, 2004; French original, 1985).

\item \verb|[12]| Jean Petitot, ``A Morphodynamical Schematization of
the Canonical Formula for Myths,'' in \textit{ The Double Twist:
From Ethnography to Morphodynamics }, edited by Pierre Maranda
(University of Toronto Press, Toronto and Buffalo, 2001, pp.
267-311).

\item \verb|[13]| Gilles Deleuze and F\'elix Guatteri, \textit{ A
Thousand Plateaus:  Capitalism and Schizophrenia} (University of
Minnesota Press, Minneapolis and London, 1987; French original,
1980).

\item \verb|[14]| Robert P. C. de Marrais, ``Flying Higher Than A
Box-Kite,'' \newline arXiv:math.RA/0207003

\item \verb|[15]| Robert P. C. de Marrais, ``Placeholder Substructures II:
Meta-Fractals, Made of Box-Kites, Fill Infinite-Dimensional Skies,''
arXiv:0704.0026 \newline[math.RA]

\item \verb|[16]| Robert P. C. de Marrais, ``Placeholder Substructures III:
A Bit-String-Driven `Recipe Theory' for Infinite-Dimensional
Zero-Divisor Spaces,'' \newline arXiv:0704.0112 [math.RA]

\end{description}

\end{document}